\pdfoutput=1
\makeatletter
\def\input@path{{03/styles/}{styles/}}
\makeatother
\documentclass[11pt]{article}
\IfFileExists{03/prefix.tex}{
\providecommand{\tcite}[2][]{\if\relax\detokenize{#1}\relax\cite{#2}\else\cite[#1]{#2}\fi}%
\providecommand{\tciteY}[2]{\cite[#1]{#2}}%

\usepackage[T1]{fontenc}%
\IfFileExists{lmodern.sty}{\usepackage{lmodern}}{}%
\usepackage[in]{fullpage}
\usepackage{amsmath}%
\usepackage{amssymb}%
\usepackage{mathtools}%
\usepackage[table]{xcolor}%
\usepackage{microtype}%
\usepackage{mleftright}
\usepackage{xspace}%
\usepackage{graphicx}%
\graphicspath{{03/}{./}}%
\usepackage{booktabs}%
\usepackage[inline]{enumitem}%
\usepackage{titlesec}%
\titlelabel{\thetitle. }%

\usepackage[amsmath,thmmarks,hyperref]{ntheorem}%
\usepackage{pgffor}%

\usepackage{hyperref}%
\usepackage[ocgcolorlinks]{ocgx2}%

\hypersetup{%
      unicode,
      breaklinks,%
      colorlinks=true,%
      urlcolor=[rgb]{0.25,0.0,0.0},%
      linkcolor=[rgb]{0.5,0.0,0.0},%
      citecolor=[rgb]{0,0.2,0.445},%
      filecolor=[rgb]{0,0,0.4},%
      anchorcolor=[rgb]{0.0,0.1,0.2}%
}

\theoremseparator{.}%

\theoremstyle{plain}%
\theoremheaderfont{\normalfont\bfseries}%
\theorembodyfont{\itshape}%
\newtheorem{theorem}{Theorem}[section]

\foreach \env/\display in {
   lemma/Lemma, conjecture/Conjecture, corollary/Corollary, claim/Claim, invariant/Invariant,
   proposition/Proposition, prop/Proposition, openproblem/Open Problem%
}{ \edef\temp{{\noexpand\newtheorem{\env}[theorem]{\display}}} \temp }

\theoremseparator{}%

\theoremseparator{.}%

\definecolor{darkslate}{rgb}{0.2, 0.2, 0.2}

\theoremseparator{.}%
\theoremheaderfont{\sffamily\bfseries\color{darkslate}}%
\theorembodyfont{\upshape}%
\foreach \env/\display in {%
   defn/Definition, definition/Definition, example/Example, exercise/Exercise, xca/Exercise, problem/Problem,
   question/Question, assumption/Assumption%
}{\edef\temp{{\noexpand\newtheorem{\env}[theorem]{\display}}} \temp}

\theoremstyle{plain}%
\theoremseparator{.}%
\theoremheaderfont{\sffamily}%
\theorembodyfont{\upshape}%
\newtheorem*{remark:unnumbered}{Remark}%

\foreach \env/\display in {
   remark/Remark, note/Note, fact/Fact, observation/Observation%
}{ \edef\temp{{\noexpand\newtheorem{\env}[theorem]{\display}}} \temp }

\newcommand{\myqedsymbol}{\ensuremath{\blacksquare}}

\theoremheaderfont{\itshape}%
\theorembodyfont{\upshape}%
\theoremstyle{nonumberplain}%
\theoremseparator{}%
\theoremsymbol{\myqedsymbol}%
\newtheorem{proof}{Proof:}%
\newtheorem{proofof}{Proof of\!}%

\definecolor{blue25emph}{rgb}{0, 0, 0.85}
\definecolor{almostblack}{rgb}{0, 0, 0.3}

\providecommand{\emphi}[1]{}
\renewcommand{\emphi}[1]{\textcolor{blue25emph}{\textbf{\emph{#1}}}}

\providecommand{\emphw}[1]{}%
\renewcommand{\emphw}[1]{{\textcolor{almostblack}{\emph{#1}}}}%

\providecommand{\emphOnly}[1]{}%
\renewcommand{\emphOnly}[1]{\emph{\textcolor{blue25emph}{\textbf{#1}}}}

\newcommand{\DrorThanks}[1]{%
   \thanks{%
      Google Research; %
      \href{mailto:aigerd@google.com}{aigerd@google.com}. %
      #1%
   }%
}
\newcommand{\SarielThanks}[1]{%
   \thanks{%
      School of Computing and Data Science; %
      University of Illinois; %
      201 N. Goodwin Avenue; %
      Urbana, IL, 61801, USA; %
      \href{mailto:sariel@illinois.edu}{sariel@illinois.edu}; %
      \url{http://sarielhp.org/}. %
      #1%
   }%
}
\newcommand{\GoranThanks}[1]{%
   \thanks{%
      Google Research; %
      \href{mailto:goranzuzic@google.com}{goranzuzic@google.com}. %
      #1%
   }%
}

\providecommand{\HLinkY}[2]{\hyperref[#2]{#1}}
\newcommand{\HLink}[2]{\hyperref[#2]{#1~\ref*{#2}}}
\newcommand{\HLinkSuffix}[3]{\hyperref[#2]{#1\ref*{#2}{#3}}}

\newcommand{\figlab}[1]{\label{fig:#1}}
\newcommand{\figref}[1]{\HLink{Figure}{fig:#1}}

\newcommand{\thmlab}[1]{{\label{theo:#1}}}
\newcommand{\thmref}[1]{\HLink{Theorem}{theo:#1}}
\newcommand{\thmrefY}[2]{\HLinkY{#2}{theo:#1}}

\newcommand{\lemlab}[1]{\label{lemma:#1}}
\newcommand{\lemref}[1]{\HLink{Lemma}{lemma:#1}}%
\newcommand{\lemrefY}[2]{\hyperref[lemma:#1]{#2}}

\newcommand{\proplab}[1]{\label{prop:#1}}
\newcommand{\propref}[1]{\HLink{Proposition}{prop:#1}}%

\newcommand{\obslab}[1]{\label{obs:#1}}
\newcommand{\obsref}[1]{\HLink{Observation}{obs:#1}}%

\newcommand{\factlab}[1]{\label{fact:#1}}
\newcommand{\factref}[1]{\HLink{Fact}{fact:#1}}%

\newcommand{\seclab}[1]{\label{sec:#1}}
\newcommand{\secref}[1]{\HLink{Section}{sec:#1}}

\providecommand{\deflab}[1]{\label{def:#1}}
\newcommand{\defref}[1]{\HLink{Definition}{def:#1}}
\newcommand{\defrefY}[2]{\hyperref[def:#1]{#2}}

\providecommand{\eqlab}[1]{}%
\renewcommand{\eqlab}[1]{\label{equation:#1}}
\newcommand{\Eqref}[1]{\HLinkSuffix{Eq.~(}{equation:#1}{)}}

\newcommand{\pth}[1]{\mleft(#1\mright)}%
\newcommand{\brc}[1]{\left\{ {#1} \right\}}
\newcommand{\set}[1]{\brc{#1}}%
\newcommand{\Set}[2]{\left\{ #1 \;\middle\vert\; #2 \right\}}
\newcommand{\sSet}[2]{\{ #1 \mid #2 \}}

\newcommand{\R}{\mathbb{R}}%

\newlist{compactenumA}{enumerate}{5}%
\setlist[compactenumA]{itemsep=-0.5ex,topsep=0.5ex,partopsep=1ex,parsep=1ex,label=(\Alph*)}%

\newlist{compactenuma}{enumerate}{5}%
\setlist[compactenuma]{itemsep=-0.5ex,topsep=0.5ex,partopsep=1ex,parsep=1ex,label=(\alph*)}%

\newlist{compactenumI}{enumerate}{5}%
\setlist[compactenumI]{itemsep=-0.5ex,topsep=0.5ex,partopsep=1ex,parsep=1ex,label=(\Roman*)}%

\newlist{compactenumi}{enumerate}{5}%
\setlist[compactenumi]{itemsep=-0.5ex,topsep=0.5ex,partopsep=1ex,parsep=1ex,label=(\roman*)}%

\newlist{compactitem}{itemize}{5}%
\setlist[compactitem]{itemsep=-0.5ex,topsep=0.5ex,partopsep=1ex,parsep=1ex,label=\ensuremath{\bullet}}%

\newcommand{\etal}{\textit{et~al.}\xspace}

\numberwithin{figure}{section}%
\numberwithin{table}{section}%
\numberwithin{equation}{section}%

\setlength{\emergencystretch}{2em}

\providecommand{\R}{\mathbb{R}}
\providecommand{\pth}[1]{\left(#1\right)}
\providecommand{\set}[1]{\left\{#1\right\}}

\newcommand{\abs}[1]{\left\lvert#1\right\rvert}

\newcommand{\conv}{\operatorname{conv}}

\newcommand{\interior}{\operatorname{int}}

\newcommand{\tr}{\operatorname{tr}}
\newcommand{\vol}{\operatorname{vol}}
\newcommand{\Sym}{\operatorname{Sym}}
\newcommand{\GL}{\operatorname{GL}}
\newcommand{\Orth}{\operatorname{O}}
\newcommand{\Bclosed}{\overline{B}}
\newcommand{\Shad}[2]{\mathsf{sh}(#1, #2)}
\newcommand{\SegFree}[1]{\mathcal{F}(#1)}
\newcommand{\GenDir}[1]{\mathcal{G}(#1)}
\newcommand{\ind}[1]{\mathbf{1}_{#1}}
\newcommand{\iprod}[2]{\langle #1, #2 \rangle}
\newcommand{\norm}[1]{\left\lVert #1 \right\rVert}
\newcommand{\lenCY}[2]{\abs{#2}_{#1}}
\newcommand{\hsupp}[1]{h^{}_{#1}}
\newcommand{\Sphere}{S^2}
\newcommand{\CSym}{\Delta}

\DeclareMathOperator{\ExC}{\mathbb{E}}
\newcommand{\Ex}[1]{\ExC\mleft[\mkern 1.5mu #1 \mkern 1.5mu\mright]}
\NewDocumentCommand{\ExY}{O{} m m}{%
    \ExC^{}_{#2}%
    \IfBlankTF{#1}%
        {\mleft[\mkern 1.5mu #3 \mkern 1.5mu\mright]}%
        {\mathopen{#1[}\mkern 2mu #3 \mkern 2mu\mathclose{#1]}}%
}
\newcommand{\mw}[1]{\overline{w}\pth{#1}}
\newcommand{\Zuzic}{\v{Z}u\v{z}i\'c\xspace}
}{}%

\begin{document}

\title{Hadwiger's Covering Conjecture in Three Dimensions}

\author{%
    Dror Aiger\DrorThanks{}%
    \and%
    Sariel Har-Peled\SarielThanks{}%
    \and%
    Goran \Zuzic\GoranThanks{}%
}

\date{\today}

\maketitle

\begin{abstract}
    Let $K \subset \R^3$ be an arbitrary convex body. We prove that $K$
    can be covered by eight smaller positive homothetic copies
    of itself: there exist a scalar $\lambda \in (0, 1)$ and translation
    vectors $v_1, \dots, v_8 \in \R^3$ such that
    \begin{equation*}
        K \subseteq \bigcup_{i=1}^{8} (\lambda K + v_i).
    \end{equation*}
    This resolves Hadwiger's covering conjecture in three dimensions.
\end{abstract}

\section{Introduction}
\seclab{intro}

A subset $K \subset \R^d$ is a \emph{convex body} if it is compact, convex,
and has nonempty interior. Denote by $c(K)$ the smallest positive integer
$m$ for which $K$ can be covered by $m$ smaller positive homothetic copies
$\lambda_1 K + v_1, \dots, \lambda_m K + v_m$ with $\lambda_i \in (0, 1)$
and $v_i \in \R^d$. Every $d$-dimensional \emph{affine cube} (a
parallelepiped, i.e., an invertible affine image of $[0, 1]^d$) satisfies
$c(K) = 2^d$, since a smaller positive homothet of $[0, 1]^d$ can contain
at most one of its $2^d$ vertices. Hadwiger~\cite{Hadwiger57}, and
independently Gohberg and Markus~\cite{GM60}, conjectured that $c(K) \le
2^d$ holds for every convex body $K \subset \R^d$; it is further
conjectured that equality holds only for affine cubes. Boltyanski~\cite{Boltyanski60} gave an equivalent
formulation in terms of boundary illumination: $c(K)$ equals the smallest
number of nonzero directions $u_1, \dots, u_m \in \R^d \setminus \set{0}$
such that for every boundary point $x \in \partial K$, the open ray
$\Set{x + t u_i}{t > 0}$ meets $\interior(K)$ for at least one $i$.

\paragraph{Background.}
The planar case was settled even before the general conjecture was stated:
Levi~\cite{Levi55} proved that every convex body $K \subset \R^2$ satisfies
$c(K) \le 4$, with equality if and only if $K$ is a parallelogram, and
observed that every smooth convex body in $\R^d$ satisfies $c(K) \le d + 1$.
For nearly seven decades, however, Hadwiger's conjecture has remained open
in every dimension $d \ge 3$.

In general dimension $d$, combining Rogers's covering density
estimate~\cite{Rogers57} with the Rogers--Shephard inequality yields the
classical upper bound $c(K) \le \binom{2d}{d} d(\ln d + \ln \ln d + 5) =
O(4^d \sqrt{d}\ln d)$ for $d \ge 3$ (see Rogers and Zong~\cite{RZ97}).
This $4^d$ barrier stood for over sixty years until Huang
\etal~\cite{HSTV22} achieved a sub-exponential improvement
$c(K) \le e^{-\Omega(\sqrt{d})}\,4^d$ using thin-shell estimates.
Subsequently, Campos \etal~\cite{CvHMT24} linked the problem to Bourgain's
slicing conjecture, obtaining an almost-exponential improvement over $4^d$
and showing that a bounded isotropic constant yields an exponential
improvement $c(K) \le (4 - \varepsilon)^d$ for an absolute constant
$\varepsilon > 0$---which now holds unconditionally following Klartag and
Lehec's resolution of the slicing problem~\cite{KL25}. For broader surveys
of covering and illumination problems, see the books of Boltyanski
\etal~\cite{BMS97} and Brass \etal~\cite{BMP05}, and the survey of Bezdek
and Khan~\cite{BK18}.

In three dimensions ($d = 3$), where the conjectured bound is $2^3 = 8$,
progress has followed two tracks. For general convex bodies without
symmetry assumptions, Lassak proved $c(K) \le 34$~\cite{Lassak88} and
later $c(K) \le 20$~\cite{Lassak98}, Papadoperakis~\cite{Papadoperakis99}
reduced the bound to $c(K) \le 16$, and Prymak~\cite{Prymak23} improved it
to $c(K) \le 14$. On the other hand, the conjectured bound $c(K) \le 8$
was established for symmetric families: centrally symmetric bodies
(Lassak~\cite{Lassak84}) and bodies symmetric about a plane
(Bezdek~\cite{Bezdek91} for polytopes and Dekster~\cite{Dekster00} in
general). Even the stronger bound $c(K) \le 6$ was proved for several
geometric classes in $\R^3$: zonotopes and more generally belt bodies
other than affine cubes (Martini~\cite{Martini87}; Boltyanski and
Martini~\cite{BM01}), bodies of constant width
(Wei{\ss}bach~\cite{Weissbach96} and Lassak~\cite{Lassak97}, following
Schramm's exponential bounds in higher dimensions~\cite{Schramm88}), and
cap bodies (Ivanov and Strachan~\cite{IS21} in the centrally symmetric
case and Arman \etal~\cite{AKP26} in general).

\paragraph{Our result.}
In this paper, we resolve Hadwiger's covering conjecture for all convex
bodies in $\R^3$.

\begin{theorem}[Hadwiger's covering conjecture in $\R^3$]
    \thmlab{main}
    Let $K \subset \R^3$ be a compact convex set with nonempty interior.
    There exist a scalar $\lambda \in (0, 1)$ and
    translation vectors $v_1, \dots, v_8 \in \R^3$ such that
    \begin{equation*}
        K \subseteq \bigcup_{i=1}^{8} (\lambda K + v_i).
    \end{equation*}
    In particular, $c(K) \le 8$.
\end{theorem}

\paragraph{Proof overview.}
We use Boltyanski's illumination formulation
(\propref{equiv-illumination}), which is affine-invariant
(\factref{affine-inv}). For any unit direction $\xi \in \Sphere$,
the directions $\pm\xi$ immediately illuminate every boundary point of $K$
outside the \emph{shadow boundary} $\Shad{K}{\xi}$---the points where a
line parallel to $\xi$ is tangent to $K$ (\figref{illumination-scheme}~(a)).
By Ewald \etal~\cite{ELR70}, almost every direction $\xi \in \Sphere$ is
parallel to no line segment on $\partial K$ (\emph{segment-free}), in
which case $\Shad{K}{\xi}$ is a simple closed curve, leaving six
directions for $\Shad{K}{\xi}$. Measuring lengths in the norm
$\norm{\cdot}_{\CSym}$ of the central symmetral
$\CSym = \frac{1}{2}(K - K)$, we show that any boundary arc of
$\CSym$-length strictly less than $2$ fits inside a smaller positive
homothet of $K$ (even when $K$ is asymmetric) and is therefore illuminated
by one direction (\lemref{curve-covering}, \figref{illumination-scheme}~(b)).
Thus it suffices to find a segment-free direction $\xi$ whose shadow
boundary has $\CSym$-length strictly less than $12$.

\begin{figure}[!htb]
    \centering
    \includegraphics[width=0.88\textwidth]{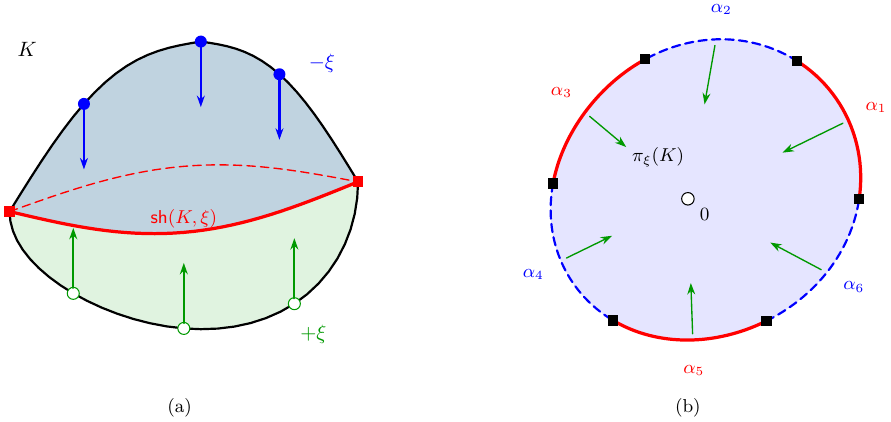}
    \caption{The $2 + 6 = 8$ illumination scheme: (a) $-\xi$ and $+\xi$
        illuminate $\partial K \setminus \Shad{K}{\xi}$, and (b)
        $\Shad{K}{\xi}$ splits into six contiguous arcs $\alpha_1, \dots,
        \alpha_6$ of $\CSym$-length strictly less than $2$, each illuminated by
    one inward direction.}
    \figlab{illumination-scheme}
\end{figure}

Affine cubes are handled directly (\obsref{parallelepipeds}), so assume
$K$ is not an affine cube. To find a short shadow boundary, we place
$\CSym$ in John's position~\cite{John48} so that $B^3 \subseteq \CSym$ is
its largest inscribed ellipsoid. Then $\norm{\cdot}_{\CSym} \le
\norm{\cdot}_2$, and by Schechtman and
Schmuckenschl{\"a}ger~\cite{SS95} and Barthe~\cite{Barthe98} its mean
width satisfies $\mw{\CSym} \le 3$, with equality only when $\CSym$ (and
hence $K$) is an affine cube; thus $\mw{\CSym} < 3$. For a polytope $P$,
combining Cauchy's perimeter formula with the average projection factor
of each edge equates the average Euclidean length of $\Shad{P}{\xi}$
over $\xi \in \Sphere$ to $4\,\mw{P}$, and central symmetrization preserves
mean width ($\mw{P} = \mw{\frac{1}{2}(P - P)}$).
We cannot simply cover approximating polytopes $P_k \to K$ and take a
limit (homothety ratios may approach $1$); instead, we use outer polytopes
$P_k \to K$ solely to transfer the average bound $4\,\mw{\CSym} < 12$ to
$\Shad{K}{\xi}$ on $K$ itself for a single segment-free direction $\xi$
(\thmref{short-shadow}).

\paragraph{AI Disclosure.}
An experimental harness at Google on top of Gemini 4 Carbon
discovered the proof with minimal human steering from Dror Aiger and
Goran \Zuzic, and an agentic harness on Gemini 4 Carbon produced the
formal verification in Lean 4. Dror Aiger, Sariel Har-Peled, and Goran
\Zuzic verified the mathematics and prepared the presentation of the
paper. In particular, real effort and time were spent on clarifying and
improving the writeup generated by the AI.

\paragraph{Organization.}
\secref{preliminaries} reviews notation, standard tools from convex
geometry, and the equivalence between covering and illumination.
\secref{proof-main} contains the proof of \thmref{main}: we establish the
curve-covering lemma in \secref{short-curves}, shadow-boundary
illumination in \secref{shadow-segments}, the mean-width bound in John's
position in \secref{john-mean-width}, and the existence of a short shadow
boundary in \secref{shadow-length}, before assembling the proof in
\secref{main-proof}. \secref{concluding} offers concluding remarks and
open questions. For completeness, \HLink{Appendix}{sec:appendix} supplies
self-contained proofs of the illumination equivalence, John's
decomposition, and the mean-width rigidity theorem.

\section{Preliminaries and illumination equivalence}
\seclab{preliminaries}

This section reviews standard notation and background tools from convex
geometry, and may be skipped by a knowledgeable reader.

For an integer $m \ge 1$, write $[m] = \set{1, \dots, m}$. Equip $\R^d$
($d \ge 1$) with the standard Euclidean inner product
$\iprod{x}{y} = \sum_{j=1}^d x_j y_j$ and norm
$\norm{x}_2 = \sqrt{\iprod{x}{x}}$, and denote the standard basis of
$\R^3$ by $e_1, e_2, e_3$.

For $x \in \R^d$ and $r > 0$, let $B(x, r)$ and $\Bclosed(x, r)$ denote the
open and closed \emphw{balls}, respectively.  In $\R^3$, write
$B^3 = \Bclosed(0, 1)$ for the closed unit ball and $\Sphere = \partial B^3$ for its
boundary which is the \emphw{unit sphere}. For a unit vector $\xi \in \Sphere$, let
$\xi^\perp = \Set{z \in \R^3}{\iprod{z}{\xi} = 0}$ and let $\pi_\xi(x) =
x - \iprod{x}{\xi}\xi$ denote orthogonal projection onto $\xi^\perp$.

For an affine subspace $W \subseteq \R^d$, a subscript $W$ on a ball, interior,
or boundary means that $W$ is the ambient space, with the induced
Euclidean metric. For $E \subseteq W$, a point $z \in E$ belongs to its relative
interior $\interior_W(E)$ if and only if $B(z,r) \cap W \subseteq E$ for some
$r > 0$. Write $\GL_3(\R)$ for the group of invertible $3 \times 3$ real
matrices and
$\Orth_3(\R) = \sSet{\smash{Q \in \R^{3 \times 3}}}{Q^\top Q = I_3}$ for
the orthogonal group, which consists of rotations and their compositions
with central inversion $x \mapsto -x$.

Let $\vol_3$ denote volume in $\R^3$, and let $\sigma$ denote surface area
on $\Sphere$ divided by $4\pi$, so that $\sigma(\Sphere) = 1$ is the uniform
probability measure on $\Sphere$. This measure is invariant under orthogonal
transformations. For any $\sigma$-integrable function $g$ on $\Sphere$, we
write $\ExY{\xi}{g(\xi)} = \int_{\Sphere} g(\xi) \, d\sigma(\xi)$ (or simply
$\Ex{g(\xi)}$ when the random direction is clear). In spherical
coordinates,
every integrable function $g$ on $\Sphere$ satisfies
\begin{equation}
    \int_{\Sphere} g(u) \, d\sigma(u)
    = \frac{1}{4\pi} \int_0^{2\pi} \int_0^\pi g(\sin\psi\cos\phi,\,
    \sin\psi\sin\phi,\, \cos\psi) \sin\psi \, d\psi \, d\phi.
    \eqlab{spherical-integral}
\end{equation}
By Archimedes' hat-box theorem, a cap of height $h$ on $\Sphere$ has
normalized area $h/2$. For $u_0 \in \Sphere$ and $\delta \in [0,2]$, the cap
within Euclidean distance $\delta$ of $u_0$ has height $\delta^2/2$, so
\begin{equation}
    \sigma\pth{\Sphere \cap \Bclosed(u_0, \delta)} = \frac{\delta^2}{4}.
    \eqlab{cap-measure}
\end{equation}
Every great circle $n^\perp \cap \Sphere$ ($n \in \Sphere$) satisfies
$\sigma(n^\perp \cap \Sphere) = 0$.

\begin{lemma}[Fatou's lemma]
    \lemlab{fatou}
    For any sequence of nonnegative measurable functions $f_1, f_2, \dots
    : \Sphere \to [0, \infty)$,
    \begin{equation*}
        \ExY[\Big]{\xi}{\liminf_{k \to \infty} f_k(\xi)}
        \le \liminf_{k \to \infty} \ExY{\xi}{f_k(\xi)}.
    \end{equation*}
\end{lemma}

\begin{definition}
    \deflab{convex-body}
    A subset $K \subset \R^d$ is a \emphi{convex body} if $K$ is compact,
    convex, and has nonempty interior $\interior(K) \neq \emptyset$.
\end{definition}

\begin{definition}
    \deflab{covering-number}
    For a convex body $K \subset \R^3$, a \emphi{smaller positive
    homothet} of $K$ is a set of the form $\lambda K + v$ with $\lambda
    \in (0, 1)$ and $v \in \R^3$. The \emphi{covering number} $c(K)$ is
    the smallest integer $m \ge 1$ for which there exist $\lambda_1,
    \dots, \lambda_m \in (0, 1)$ and $v_1, \dots, v_m \in \R^3$ such that
    $K \subseteq \bigcup_{i=1}^m (\lambda_i K + v_i)$.
\end{definition}

\begin{definition}
    \deflab{parallelepiped}
    A subset $K \subset \R^3$ is an \emphi{affine cube} (or \emph{parallelepiped}) if $K = M([-1,
    1]^3) + b$ for some invertible matrix $M \in \GL_3(\R)$ and vector $b
    \in \R^3$.
\end{definition}

\begin{definition}
    \deflab{illumination}
    Let $K \subset \R^3$ be a convex body. A nonzero vector $u \in \R^3
    \setminus \set{0}$ \emphi{illuminates} a boundary point $x \in
    \partial K$ if there exists $t > 0$ such that $x + t u \in
    \interior(K)$. A collection of nonzero vectors $u_1, \dots, u_m \in
    \R^3 \setminus \set{0}$ \emphi{illuminates $\partial K$} if every $x
    \in \partial K$ is illuminated by at least one $u_i$.
\end{definition}

The key property of illumination of the boundary is that it requires
every point to be illuminated from ``above'', and not from the side. That
is, a point that is on the silhouette of the convex body when projected
in a direction $u$ is not illuminated by this direction (\figref{illumination-equiv}~(a)).

\begin{figure}[!htb]
    \centering
    \includegraphics[width=0.88\textwidth]{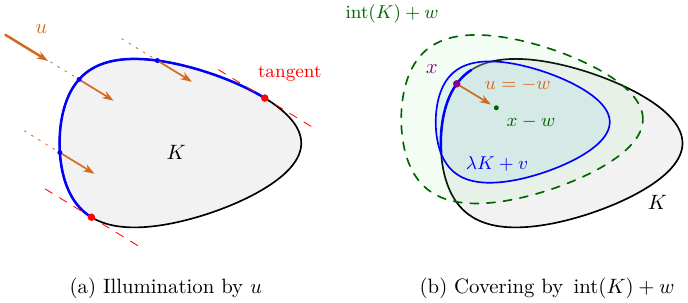}
    \caption{Boundary illumination and covering equivalence
        (\propref{equiv-illumination}): (a) $u$ illuminates points whose
        inward ray meets $\interior(K)$, missing the tangent points; (b)
        $\lambda K + v \subseteq \interior(K) + w$ implies $u = -w$
    illuminates the covered arc.}
    \figlab{illumination-equiv}
\end{figure}

\begin{definition}
    \deflab{support-function}
    For a nonempty compact subset $C \subset \R^d$, its \emphi{support
    function} $\hsupp{C} : \R^d \to \R$ is defined by $\hsupp{C}(u) = \max_{x \in C}
    \iprod{x}{u}$. In $\R^3$, the width of $C$ in direction $u \in \Sphere$
    is $\hsupp{C}(u) + \hsupp{C}(-u)$, and its \emphi{mean width} is
    \begin{equation*}
        \mw{C} = \ExY{u}{\hsupp{C}(u) + \hsupp{C}(-u)} = 2\,\ExY{u}{\hsupp{C}(u)},
    \end{equation*}
    where the second equality uses the symmetry of $\sigma$ under $u \mapsto -u$.
\end{definition}

\begin{definition}
    \deflab{minkowski-gauge}
    For a convex body $C \subset \R^d$ satisfying $0 \in \interior(C)$,
    its \emphi{Minkowski gauge} $\mu_C : \R^d \to [0, \infty)$ is defined
    by $\mu_C(x) = \inf\Set{t > 0}{x \in t C}$. When $C = -C$ is
    centrally symmetric, we also write $\norm{x}_C = \mu_C(x)$ and call
    $\norm{\cdot}_C$ the \emphi{Minkowski norm} of $C$.
\end{definition}

\begin{lemma}[{Properties of the Minkowski gauge~\cite{Schneider14}}]
    \lemlab{gauge-props}
    Let $C \subset \R^d$ be a convex body with $0 \in \interior(C)$, and
    choose $0 < r \le R$ such that $\Bclosed(0, r) \subseteq C \subseteq
    \Bclosed(0, R)$. For all $x, y \in \R^d$ and $\alpha \ge 0$, we have
    $\frac{1}{R}\norm{x}_2 \le \mu_C(x) \le \frac{1}{r}\norm{x}_2$,
    $\mu_C(\alpha x) = \alpha \mu_C(x)$, $\mu_C(x + y) \le \mu_C(x) +
    \mu_C(y)$, $\abs{\mu_C(x) - \mu_C(y)} \le \frac{1}{r}\norm{x
    - y}_2$, $x \in \lambda C \iff \mu_C(x) \le \lambda$ for all
    $\lambda > 0$, and $x \in \interior(C) \iff \mu_C(x) < 1$. In
    particular, if $C = -C$, then $\norm{\cdot}_C = \mu_C$ is a norm on
    $\R^d$.
\end{lemma}

\begin{definition}
    \deflab{difference-body}
    For a convex body $K \subset \R^3$, its \emphi{central symmetral}
    (or \emph{halved difference body}) is the set $\CSym = \frac{1}{2}(K
    - K) = \sSet{\frac{1}{2}(x - y)}{x, y \in K}$. Exchanging $x$ and $y$
    shows that $\CSym = -\CSym$ (\figref{halved-difference-body});
    moreover, $\CSym$ is a centrally symmetric convex body with $0 \in
    \interior(\CSym)$, so $\norm{\cdot}_{\CSym}$ is a norm on $\R^3$
    with unit ball $\CSym$. Its support function is
    $\hsupp{\CSym}(u) = \frac{1}{2}(\hsupp{K}(u) + \hsupp{K}(-u))$, and
    hence $\mw{\CSym} = \mw{K}$.
\end{definition}

\begin{figure}[!htb]
    \centering
    \includegraphics[width=0.75\textwidth]{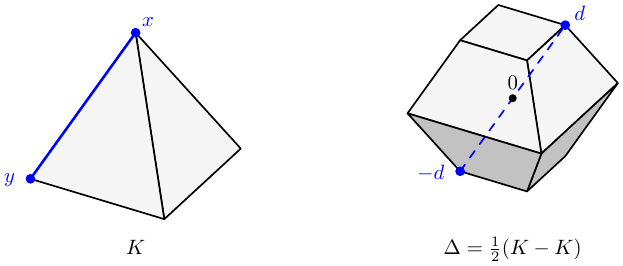}
    \caption{A square pyramid $K$ and its central symmetral $\CSym =
        \frac{1}{2}(K - K)$, a truncated square bipyramid. The
        highlighted slant edge of $K$ gives the opposite vertices
        $d=(x-y)/2$ and $-d$ of $\CSym$.}
    \figlab{halved-difference-body}
\end{figure}

\begin{definition}
    \deflab{curve-variation}
    Let $C \subset \R^3$ be a centrally symmetric convex body with $0 \in
    \interior(C)$, and let $\gamma : [\alpha, \beta] \to \R^3$ ($\alpha < \beta$) be a
    continuous curve. For $\alpha \le t_1 < t_2 \le \beta$, write $\gamma[t_1, t_2]$
    for the restricted subcurve $\gamma|_{[t_1, t_2]}$. The \emphi{$C$-length}
    of $\gamma$ on $[\alpha, \beta]$ is
    \begin{equation*}
        \lenCY{C}{\gamma}
        =
        \sup_{\alpha = t_0 < t_1 < \dots < t_N = \beta}
        \sum_{r=1}^N \norm{\gamma(t_r) - \gamma(t_{r-1})}_C \in [0, \infty],
    \end{equation*}
    and $\lenCY{C}{\gamma[t_1, t_2]}$ denotes the $C$-length of the subcurve
    $\gamma[t_1, t_2]$.
\end{definition}

\begin{definition}
    \deflab{shadow-boundary}
    For a convex body $K \subset \R^3$ and a unit vector $\xi \in \Sphere$, a line
    parallel to $\xi$ is \emphi{tangent} to $K$ if it intersects $K$ but misses
    $\interior(K)$. The \emphi{shadow boundary} of $K$ in direction $\xi$ is the
    set of all boundary points lying on a line parallel to $\xi$ that is tangent
    to $K$:
    \begin{equation*}
        \Shad{K}{\xi} = \Set{x \in \partial K}{(x + \R \xi) \cap \interior(K) = \emptyset}.
    \end{equation*}
\end{definition}

\begin{definition}
    \deflab{segment-free}
    For a convex body $K \subset \R^3$, a unit vector $\xi \in \Sphere$ is a
    \emphi{segment-free direction} of $K$ if $\partial K$ contains no
    nondegenerate line segment parallel to $\xi$, meaning there do not
    exist distinct points $x, y \in \partial K$ with $[x, y] \subseteq
    \partial K$ and $y - x \in \R \xi$. Denote by $\SegFree{K} \subseteq
    \Sphere$ the set of all segment-free directions of $K$.
\end{definition}

Intuitively, almost every direction is segment-free, though proving this
is surprisingly subtle.

\begin{theorem}[{Ewald \etal~\cite{ELR70}}]
    \thmlab{segment-free-measure}
    For every convex body $K \subset \R^3$, the set $\Sphere \setminus \SegFree{K}$
    of directions parallel to a boundary segment of $K$ has measure $0$.
\end{theorem}

\begin{lemma}[\tciteY{Section~1.3}{Schneider14}]
    \lemlab{separation}
    Let $V$ be a finite-dimensional real inner-product space equipped
    with $\norm{v} = \sqrt{\iprod{v}{v}}$.
    \begin{compactenumi}
    \item Let $Q \subseteq V$ be nonempty, closed, and convex, and let
        $y_0 \in V \setminus Q$. There exists a point $q_* \in Q$
        minimizing $\norm{y_0 - q}$ over $q \in Q$, and the nonzero
        vector $w = y_0 - q_*$ satisfies $\iprod{w}{q} \le
        \iprod{w}{y_0} - \norm{w}^2$ for all $q \in Q$.
    \item Let $C \subset \R^d$ be a convex body and let $x_0 \in \partial
        C$. There exists $n \in \R^d$ with $\norm{n}_2 = 1$ such that
        $\iprod{n}{x} \le \iprod{n}{x_0} = \hsupp{C}(n)$ for all $x \in C$.
    \end{compactenumi}
\end{lemma}

\begin{proposition}[{Equivalence of covering and illumination~\cite{Boltyanski60}}]
    \proplab{equiv-illumination}
    Let $K \subset \R^3$ be a convex body and let $m \ge 1$ be an
    integer. The following four statements are equivalent:
    \begin{compactenumi}
    \item There exist $\lambda_1, \dots, \lambda_m \in (0, 1)$ and $v_1,
        \dots, v_m \in \R^3$ such that $K \subseteq \bigcup_{i=1}^m
        (\lambda_i K + v_i)$.
    \item There exist $\lambda \in (0, 1)$ and $v_1, \dots, v_m \in \R^3$
        such that $K \subseteq \bigcup_{i=1}^m (\lambda K + v_i)$.
    \item There exist $w_1, \dots, w_m \in \R^3$ such that $K \subseteq
        \bigcup_{i=1}^m (\interior(K) + w_i)$.
    \item There exist nonzero vectors $u_1, \dots, u_m \in \R^3 \setminus
        \set{0}$ that illuminate $\partial K$.
    \end{compactenumi}
\end{proposition}

For completeness, a proof of \propref{equiv-illumination} (see also
\figref{illumination-equiv}~(b)) is given in
\HLink{Appendix}{sec:proof-equiv-illumination}.

\begin{fact}
    \factlab{affine-inv}
    For every invertible affine map $\Phi(x) = M x + b$ ($M \in
    \GL_3(\R)$ and $b \in \R^3$) and every convex body $K \subset \R^3$,
    we have $c(\Phi(K)) = c(K)$, where $c(\cdot)$ is the covering
    number (see \defref{covering-number}).
\end{fact}

\begin{observation}
    \obslab{parallelepipeds}
    Every affine cube $K \subset \R^3$ satisfies $c(K) = 8$: by
    \factref{affine-inv} we may assume $K = [0, 1]^3$, which is covered by
    eight translates of $\frac{2}{3}[0, 1]^3$, while each of its eight
    corners requires its own smaller homothet.
\end{observation}

\section{Proof of \texorpdfstring{\thmref{main}}{Theorem 1.1}}
\seclab{proof-main}

Throughout \HLink{Subsections}{sec:short-curves}--\ref{sec:shadow-length},
let $K \subset \R^3$ be a convex body with $0 \in \interior(K)$ and
central symmetral $\CSym = \frac{1}{2}(K - K)$.

\subsection{Covering short curves by a single homothet}
\seclab{short-curves}

We now show that any continuous curve of $\CSym$-length strictly less than $2$
is covered by a single smaller positive homothet of $K$, using Minkowski
addition of segment covers and compactness.

\begin{figure}[!htb]
    \centering
    \includegraphics[width=0.88\textwidth]{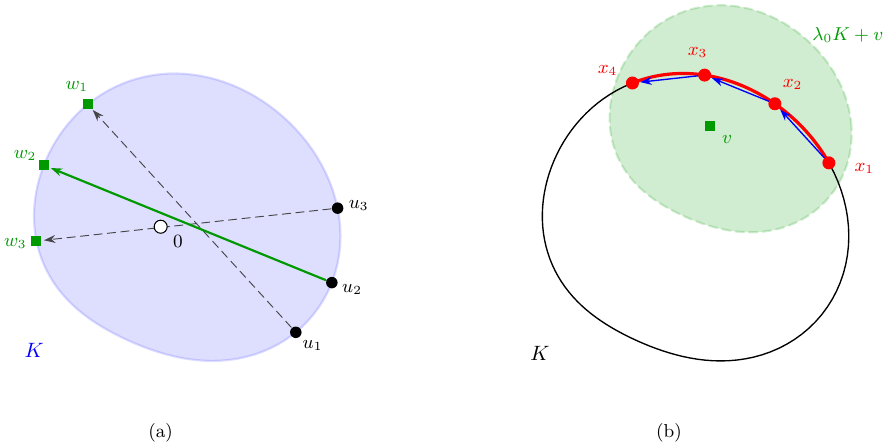}
    \caption{Covering a short arc (\lemref{curve-covering}): (a) chord
        differences $x_{r+1} - x_r = \zeta_r(w_r - u_r)$ with $u_r, w_r \in K$,
        and (b) sample points $x_1, x_2, x_3, x_4$ along $\gamma[\alpha,
    \beta]$ enclosed by a single smaller homothet $\lambda_0 K + v$.}
    \figlab{curve-covering}
\end{figure}

\begin{lemma}[Covering short curves by a single homothet]
    \lemlab{curve-covering}
    Let $\gamma : [\alpha, \beta]
    \to \R^3$ ($\alpha < \beta$) be a continuous curve satisfying
    $\lenCY{\CSym}{\gamma} < 2$, and set $\lambda_0 =
    \max\pth{\frac{1}{2}\lenCY{\CSym}{\gamma}, \frac{1}{2}} \in
    [\frac{1}{2}, 1)$. Then there exists $v \in \R^3$ such that
    $\gamma[\alpha, \beta] \subseteq \lambda_0 K + v$. If moreover $\gamma[\alpha,
    \beta] \subseteq \partial K$, then $v \neq 0$ and the vector $u = -v
    \in \R^3 \setminus \set{0}$ illuminates every point of $\gamma[\alpha,
    \beta]$ in $K$.
\end{lemma}

\begin{proof}
    Fix an integer $N \ge 1$ and partition $[\alpha, \beta]$ into $N$ equal
    subintervals with endpoints $\alpha = t_1 < t_2 < \dots < t_{N+1} =
    \beta$. Set $x_j = \gamma(t_j)$ for $j \in [N+1]$ and $d_r = x_{r+1} -
    x_r$ for $r \in [N]$ (\figref{curve-covering}). For each $r \in [N]$,
    let $u_r w_r$ be a longest chord of $K$ parallel to $x_r x_{r+1}$ and
    oriented in the same direction (taking $u_r = w_r = 0$ if $d_r = 0$).
    When $d_r \neq 0$, maximality of $u_r w_r$ places $\frac{1}{2}(w_r -
    u_r) \in \partial \CSym$, so $\norm{w_r - u_r}_{\CSym} = 2$; setting $\zeta_r =
    \norm{d_r}_{\CSym} / \norm{w_r - u_r}_{\CSym} = \frac{1}{2}\norm{d_r}_{\CSym}$ (and
    $\zeta_r = 0$ if $d_r = 0$) gives $d_r = \zeta_r(w_r - u_r)$ and
    \[
        \zeta
        =
        \sum_{r=1}^N \zeta_r
        \leq
        \frac{1}{2}\lenCY{\CSym}{\gamma} \le
        \lambda_0.
   \]

    Set $v_N = x_1 - \sum_{r=1}^N \zeta_r u_r \in x_1 - \zeta K \subseteq
    x_1 - \lambda_0 K$. For every $j \in [N+1]$, since 
    \begin{equation*}
        x_j = x_1 + \sum_{r=1}^{j-1} d_r = x_1 + \sum_{r=1}^{j-1} \zeta_r(w_r -
        u_r),
    \end{equation*}
    we have
    \begin{equation*}
        x_j - v_N
        \;=\; \sum_{r=1}^{j-1} \zeta_r w_r + \sum_{r=j}^N \zeta_r u_r
        \;\in\; \zeta K
        \;\subseteq\; \lambda_0 K
    \end{equation*}
    by convexity of $K$ (using $u_r, w_r, 0 \in K$ and
    $\zeta \le \lambda_0$). Thus the polygonal approximation
    $x_1 x_2 \dots x_{N+1}$ of $\gamma$ is contained in
    $\lambda_0 K + v_N$. Since $v_N \in x_1 - \lambda_0 K$ lies in a compact set,
    passing to a convergent subsequence $v_{N_k} \to v$ as $k \to \infty$ yields
    $\gamma[\alpha, \beta] \subseteq \lambda_0 K + v$ by continuity of $\gamma$ and as $K$ is closed.

    Finally, if $\gamma[\alpha, \beta] \subseteq \partial K$, then
    $\gamma[\alpha, \beta] - v \subseteq \lambda_0 K \subseteq \interior(K)$
    (since $0 \in \interior(K)$ and $\lambda_0 < 1$). Because
    $\partial K \cap \interior(K) = \emptyset$, this forces $v \neq 0$, and
    $u = -v$ illuminates $\gamma[\alpha, \beta]$.
\end{proof}

\subsection{Illuminating outside the shadow boundary}
\seclab{shadow-segments}

We now show that $\set{+\xi, -\xi}$ illuminates $\partial K \setminus
\Shad{K}{\xi}$ for every unit direction $\xi \in \Sphere$, and that for every
$\xi \in \SegFree{K}$ the shadow boundary $\Shad{K}{\xi}$ is a continuous
closed curve.

\begin{observation}
    \obslab{outside-shadow-obs}
    For any $\xi \in \Sphere$, the planar
    projection $D = \pi_\xi(K) \subset \xi^\perp$ is a $2$-dimensional convex
    body satisfying $\pi_\xi(\interior(K)) = \interior_{\xi^\perp}(D)$ and
    $\Shad{K}{\xi} = \partial K \cap \pi_\xi^{-1}(\partial_{\xi^\perp}(D))$. In
    particular, every point $x \in \partial K \setminus \Shad{K}{\xi}$ lies on
    a line $x + \R \xi$ meeting $\interior(K)$, and is therefore illuminated
    by $+\xi$ or $-\xi$. Moreover, if $\xi \in \SegFree{K}$, then for every
    $z \in \partial_{\xi^\perp}(D)$ the intersection $K \cap \pi_\xi^{-1}(z)$
    consists of a single point lying in $\Shad{K}{\xi}$.
\end{observation}

\begin{observation}
    \obslab{outside-shadow}
    For any $\xi \in \SegFree{K}$ and $D = \pi_\xi(K) \subset \xi^\perp$, since
    $\Shad{K}{\xi} = \partial K \cap \pi_\xi^{-1}(\partial_{\xi^\perp}(D))$ is
    compact, the continuous bijection $\pi_\xi|_{\Shad{K}{\xi}} : \Shad{K}{\xi}
    \to \partial_{\xi^\perp}(D)$ from \obsref{outside-shadow-obs} is a
    homeomorphism; in particular, $\Shad{K}{\xi}$ is a simple closed curve
    (\figref{segment-free-shadow}).
    Because $0 \in \interior_{\xi^\perp}(D)$, fixing polar coordinates in
    $\xi^\perp$ yields a continuous polar parameterization $\gamma_\xi : [0,
    2\pi] \to \Shad{K}{\xi}$ (bijective on $[0, 2\pi)$, with $\gamma_\xi(0) =
    \gamma_\xi(2\pi)$), where $\gamma_\xi(\theta)$ is the unique point of
    $\Shad{K}{\xi}$ whose projection $\pi_\xi(\gamma_\xi(\theta)) \in
    \partial_{\xi^\perp}(D)$ lies on the ray from $0$ at angle $\theta$.
    For any centrally symmetric convex body $C \subset \R^3$ with $0 \in
    \interior(C)$, the $C$-length of $\gamma_\xi$ is independent of the
    choice of simple closed parameterization of $\Shad{K}{\xi}$, so we
    write $\lenCY{C}{\Shad{K}{\xi}} = \lenCY{C}{\gamma_\xi}$.
\end{observation}

\begin{figure}[!htb]
    \centering
    \includegraphics[width=0.88\textwidth]{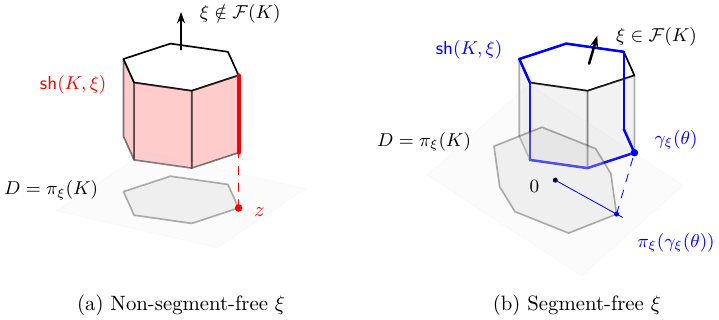}
    \caption{Role of segment-free directions (\obsref{outside-shadow}): (a)
        for $\xi \notin \SegFree{K}$, a boundary point $z \in
        \partial_{\xi^\perp}(D)$ can lift to a vertical segment; (b) for
        $\xi \in \SegFree{K}$, $\Shad{K}{\xi}$ projects homeomorphically onto
    $\partial_{\xi^\perp}(D)$.}
    \figlab{segment-free-shadow}
\end{figure}

\subsection{Mean width in John's position}
\seclab{john-mean-width}

By \factref{affine-inv}, we may apply an invertible linear transformation
placing the centrally symmetric central symmetral $\CSym = \frac{1}{2}(K
- K)$ in \emph{John's position}, meaning that $B^3 \subseteq \CSym$ is the
maximum-volume origin-centered ellipsoid inscribed in $\CSym$ (equivalently,
$\abs{\det M} \le 1$ whenever $M \in \R^{3 \times 3}$ satisfies $M(B^3)
\subseteq \CSym$). Such a position exists because $0 \in \interior(\CSym)$
and $\CSym$ is compact, so $\Set{M \in \R^{3 \times 3}}{M(B^3) \subseteq
\CSym}$ is compact and contains a positive multiple of $I_3$, whence
$\abs{\det M}$ attains a positive maximum at some $M_* \in \GL_3(\R)$
and applying $M_*^{-1}$ places $\CSym$ in John's position. We now
establish the sharp mean-width inequality and its cube rigidity case.

\begin{lemma}[\tciteY{Section~3.2}{Schneider14}]
    \lemlab{diff-parallelepiped}
    The central symmetral $\CSym$ is an affine cube if and only if $K$ is an affine cube.
\end{lemma}

\begin{proof}
    If $K = M([-1, 1]^3) + b$ for some $M \in \GL_3(\R)$ and $b \in \R^3$, then $\CSym = \frac{1}{2}(K - K) = M([-1, 1]^3)$ is an affine cube. Conversely, suppose $\CSym$ is an affine cube; applying an invertible affine map, we may assume $\CSym = [-1, 1]^3$, so $K - K = [-2, 2]^3$. Translating $K$ so that $\min_{x \in K} x_j = -1$ for each coordinate $j \in [3]$, the bound $x_j - y_j \le 2$ for all $x, y \in K$ gives $\max_{x \in K} x_j \le 1$, hence $K \subseteq [-1, 1]^3$. For every corner $\epsilon \in \set{-1, 1}^3$, we have $2\epsilon \in K - K$, so $2\epsilon = x - y$ for some $x, y \in K \subseteq [-1, 1]^3$, which forces $x = \epsilon \in K$. By convexity of $K$, $K = [-1, 1]^3$.
\end{proof}

\begin{remark:unnumbered}
    Recall that for a column vector $v \in \R^3$, the \emph{outer product}
    $v v^\top \in \R^{3 \times 3}$ is the rank-$1$ symmetric matrix with
    $(i, j)$-entry $v_i v_j$. Applied to any $x \in \R^3$, it yields
    $(v v^\top)x = v(v^\top x) = \iprod{v}{x} v$, which for a unit vector
    $v \in \Sphere$ is the orthogonal projection of $x$ onto the line spanned
    by $v$. Thus the matrix identity $\sum_{s=1}^r c_s v_s v_s^\top = I_3$
    below states that every $x \in \R^3$ decomposes as $x = \sum_{s=1}^r
    c_s \iprod{v_s}{x} v_s$.
\end{remark:unnumbered}

\begin{theorem}[John's decomposition, \tcite{John48}]
    \thmlab{john}
    Suppose $\CSym \subset \R^3$ is a centrally symmetric convex body ($\CSym =
    -\CSym$ and $0 \in \interior(\CSym)$) in John's position. Then there exist an
    integer $r \ge 3$, pairwise noncollinear contact vectors $v_1, \dots,
    v_r \in \Sphere \cap \partial \CSym$ ($v_s \neq \pm v_t$ for $s \neq t$), and
    weights $c_1, \dots, c_r > 0$ such that:
    \begin{compactenumi}
    \item $\sum_{s=1}^r c_s v_s v_s^\top = I_3$ and $\sum_{s=1}^r c_s =
        3$, and
    \item $\abs{\iprod{v_s}{x}} \le 1$ for all $x \in \CSym$ and all $s \in
        [r]$.
    \end{compactenumi}
\end{theorem}

For completeness, a proof of \thmref{john} is given in
\HLink{Appendix}{sec:proof-john}.  Geometrically, \thmref{john}~(i)--(ii)
trap $\CSym$ inside a small enclosing zonotope: every $x \in \CSym$ decomposes as
$x = \sum_{s=1}^r c_s \iprod{v_s}{x} v_s$ with
$\abs{\iprod{v_s}{x}} \le 1$, so
\begin{equation*}
    \CSym \subseteq Z = \sum_{s=1}^r [-c_s, c_s] v_s,
\end{equation*}
a Minkowski sum of $r$ segments whose half-lengths sum to $\sum_{s=1}^r
c_s = 3$ (the same total half-length, and hence the same mean width, as
the unit cube $[-1, 1]^3$). When $r = 3$, the contact vectors are
orthonormal and $Z$ is a rotated cube; when $r > 3$, the nonorthogonal
contacts leave a strict gap between $\CSym$ and $Z$ (\figref{john-cube}),
forcing the mean width of $\CSym$ strictly below that of $[-1, 1]^3$.

\begin{figure}[!htb]
    \centering
    \includegraphics[width=0.88\textwidth]{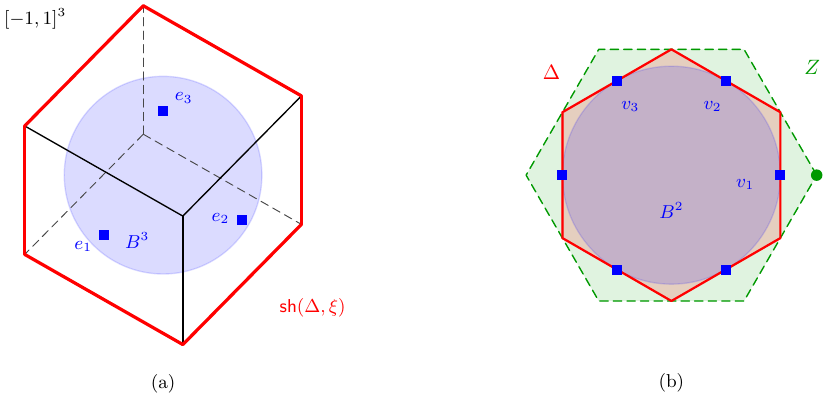}
    \caption{Mean width in John's position
        (\thmref{mean-width-rigidity}): (a) for $[-1, 1]^3$, generic shadow
        boundaries $\Shad{\CSym}{\xi}$ consist of six edges of $\CSym$-length $2$,
        whereas (b) in a $2$-dimensional analogue around $B^2$, nonorthogonal contacts
        $v_1, v_2, v_3$ leave a strict gap between $\CSym$ and the bounding
        zonotope $Z = \sum_{s=1}^r [-c_s, c_s]v_s$.}
    \figlab{john-cube}
\end{figure}

\begin{theorem}[{Mean width in John's position~\cite{SS95,Barthe98}}]
    \thmlab{mean-width-rigidity}
    When $\CSym$ is in John's position, $\mw{\CSym} \le 3$,
    with equality if and only if $\CSym = Q^\top([-1, 1]^3)$ for some $Q \in
    \Orth_3(\R)$, which holds if and only if $K$ is an affine cube.
\end{theorem}

For completeness, a proof of \thmref{mean-width-rigidity} is given in
\HLink{Appendix}{sec:proof-mean-width}.

\subsection{Average length of the shadow boundary}
\seclab{shadow-length}

We finally relate the average shadow-boundary length of an outer polytope
to its mean width via Cauchy's planar perimeter formula, pass to a short
shadow boundary on $K$ itself, and complete the proof of \thmref{main}.

\begin{lemma}[Outer polytope approximation, \tcite{Schneider14}]
    \lemlab{outer-polytope}
    For every $\varepsilon > 0$, there exists a convex polytope $P \subset
    \R^3$ with $0 \in \interior(P)$ satisfying $K \subseteq P \subseteq (1
    + \varepsilon)K$.
\end{lemma}

\begin{lemma}[Great-circle averages on the sphere]
    \lemlab{great-circle}
    Let $f : \Sphere \to \R$ be continuous. For each $\xi \in \Sphere$, the
    uniform average of $f$ over the great circle $\xi^\perp \cap \Sphere$ is
    \[
        \ExY{u \in \xi^\perp \cap \Sphere}{f(u)}
        = \frac{1}{2\pi}\int_0^{2\pi} f(\cos\alpha \, e_1(\xi) +
        \sin\alpha \, e_2(\xi)) \, d\alpha,
    \]
    which is independent of the choice of orthonormal basis $(e_1(\xi),
    e_2(\xi))$ of $\xi^\perp$ and Borel measurable in $\xi$, and satisfies
    \begin{equation*}
        \ExY[\big]{\xi}{\ExY{u \in \xi^\perp \cap \Sphere}{f(u)}}
        = \ExY{u}{f(u)}.
    \end{equation*}
\end{lemma}

This follows immediately from the $\Orth_3(\R)$-invariance of $\sigma$ and
Fubini's theorem.

\begin{observation}
    \obslab{polyhedral-silhouette}
    Let $P \subset \R^3$ be a convex polytope with $0 \in \interior(P)$,
    facets $F_1, \dots, F_M$ having outward unit normals $n_1, \dots,
    n_M \in \Sphere$, and edge set $E(P)$, where each edge $e = [p_e, q_e]
    \in E(P)$ is shared by two facets with outward unit normals $n_e,
    n_e' \in \set{n_1, \dots, n_M}$ and \emphw{exterior dihedral angle}
    \[
        \theta(e) = \arccos(\iprod{n_e}{n_e'}) \in (0, \pi).
    \]
     Define
    \begin{equation*}
        \GenDir{P}
        = \Sphere \setminus \bigcup_{i=1}^M (n_i^\perp \cap \Sphere)
        = \Set{\xi \in \Sphere}{\iprod{n_i}{\xi} \neq 0 \ \forall i \in [M]}.
    \end{equation*}
    Then $\GenDir{P}$ is open in $\Sphere$ with $\sigma(\GenDir{P}) = 1$ and
    $\GenDir{P} \subseteq \SegFree{P}$ (every segment on $\partial P$
    lies in some facet $F_i$ and is orthogonal to $n_i$). Moreover, for
    every $\xi \in \GenDir{P}$, no facet is parallel to $\xi$, so the shadow
    boundary $\Shad{P}{\xi}$ is simply a closed cycle of edges on
    $\partial P$---namely, the silhouette edges separating front-facing
    ($\iprod{n_i}{\xi} > 0$) and back-facing ($\iprod{n_i}{\xi} < 0$)
    facets:
    \begin{equation*}
        \Shad{P}{\xi} = \bigcup_{e \in E(P) : \iprod{n_e}{\xi}\iprod{n_e'}{\xi} <
        0} e.
    \end{equation*}
    Consequently, for any centrally symmetric convex body $C \subset \R^3$
    with $0 \in \interior(C)$, the $C$-length of $\Shad{P}{\xi}$
    (\obsref{outside-shadow}) is simply the sum of the $C$-lengths of
    these edges:
    \begin{equation*}
        \lenCY{C}{\Shad{P}{\xi}} = \sum_{e \in E(P)} \norm{q_e - p_e}_C \,
        \ind{\set{\iprod{n_e}{\xi}\iprod{n_e'}{\xi} < 0}}.
    \end{equation*}
\end{observation}

\begin{lemma}[Average length of a polyhedral shadow boundary]
    \lemlab{polyhedral-crofton}
    Using the notation of \obsref{polyhedral-silhouette}, for any centrally
    symmetric convex body $C \subset \R^3$ with $0 \in \interior(C)$, we
    have
    \begin{align*}
        \ExY{\xi}{\lenCY{C}{\Shad{P}{\xi}}}
        &= \frac{1}{\pi}\sum_{e \in E(P)} \norm{q_e - p_e}_C\,\theta(e) \\
        \qquad\text{and}\qquad
        \frac{1}{\pi}\sum_{e \in E(P)} \norm{q_e - p_e}_2 \theta(e)
        &= 4\,\mw{P}.
    \end{align*}
\end{lemma}

\begin{figure}[!htb]
    \centering
    \includegraphics[width=0.88\textwidth]{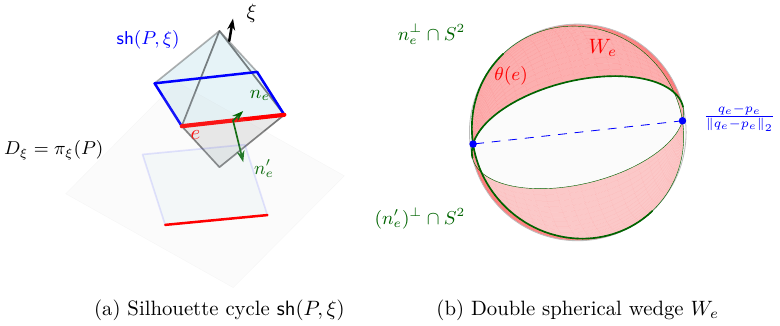}
    \caption{Polyhedral shadow-boundary average (\lemref{polyhedral-crofton}):
        (a) silhouette cycle $\Shad{P}{\xi}$ separating front and back facets
        with edge $e$ and normals $n_e, n_e'$; (b) double spherical wedge
    $W_e \subset \Sphere$ of measure $\theta(e)/\pi$.}
    \figlab{polyhedral-crofton}
\end{figure}

\begin{proof}
    For each edge $e = [p_e, q_e] \in E(P)$, the set $W_e = \sSet{\xi \in
    \Sphere}{\iprod{n_e}{\xi}\iprod{n_e'}{\xi} < 0}$ of directions in which $e$
    appears on the silhouette $\Shad{P}{\xi}$ (\figref{polyhedral-crofton}~(a))
    is the double spherical wedge
    (two antipodal lunes) bounded by the great circles $n_e^\perp \cap
    \Sphere$ and $(n_e')^\perp \cap \Sphere$, which intersect along the edge axis
    $\pm (q_e - p_e)/\norm{q_e - p_e}_2$ at angle $\theta(e)$
    (\figref{polyhedral-crofton}~(b)). Since the
    two wedges span azimuthal angle $2\theta(e)$ out of $2\pi$, we have
    $\sigma(W_e) = \theta(e)/\pi$. Taking the expectation of
    \obsref{polyhedral-silhouette} over $\xi \in \Sphere$ (since
    $\sigma(\GenDir{P}) = 1$) immediately yields the first identity.
    
    For the second identity, observe first that for any vector $v \in
    \R^3$, the spherical average of its projected length onto $\xi^\perp$ is
    \[
        \ExY{\xi}{\norm{\pi_\xi(v)}_2} = \norm{v}_2 \cdot
        \frac{1}{2}\int_0^\pi \sin^2\psi \, d\psi = \frac{\pi}{4}\norm{v}_2
    \]
    (by \Eqref{spherical-integral}, measuring $\psi$ from $v$). By
    rotational symmetry around the edge axis $q_e - p_e$ of the double
    wedge $W_e$, restricting this expectation to $W_e$ multiplies it by
    $\sigma(W_e) = \theta(e)/\pi$, yielding
    \begin{equation*}
        \ExY{\xi}{\norm{\pi_\xi(q_e - p_e)}_2\,\ind{W_e}(\xi)}
        = \frac{\theta(e)}{\pi} \cdot \frac{\pi}{4}\norm{q_e - p_e}_2
        = \frac{1}{4}\norm{q_e - p_e}_2 \theta(e).
    \end{equation*}
    On the other hand, for every $\xi \in \GenDir{P}$, the boundary of the
    planar convex polygon $D_\xi = \pi_\xi(P)$ is the projection of
    $\Shad{P}{\xi}$, so its perimeter is
    $\sum_{e \in E_P(\xi)} \norm{\pi_\xi(q_e - p_e)}_2$, where
    $E_P(\xi) = \sSet{e \in E(P)}{\xi \in W_e}$. By Cauchy's planar perimeter
    formula~\cite{Schneider14} (with
    $u(\alpha) = \cos\alpha \, e_1(\xi) + \sin\alpha \, e_2(\xi)$ and
    $\hsupp{D_\xi} = \hsupp{P}$ on $\xi^\perp$),
    \begin{align*}
        \sum_{e \in E_P(\xi)} \norm{\pi_\xi(q_e - p_e)}_2
        &= \int_0^\pi \pth{\Bigl.\hsupp{D_\xi}(u(\alpha)) +
        \hsupp{D_\xi}(-u(\alpha))} \, d\alpha \\
        &= \int_0^{2\pi} \hsupp{P}(u(\alpha)) \, d\alpha
        = 2\pi\,\ExY{u \in \xi^\perp \cap \Sphere}{\hsupp{P}(u)}.
    \end{align*}
    Taking the expectation over $\xi \in \GenDir{P}$---using
    $\sigma(\GenDir{P}) = 1$, $e \in E_P(\xi) \iff \xi \in W_e$, and
    \lemref{great-circle} for $f = \hsupp{P}$---we
    obtain
    \begin{align*}
        \frac{1}{4}\sum_{e \in E(P)} \norm{q_e - p_e}_2 \theta(e)
        &= \sum_{e \in E(P)} \ExY[\big]{\xi}{\norm{\pi_\xi(q_e - p_e)}_2\,\ind{W_e}(\xi)} \\
        &= \ExY[\bigg]{\xi}{\sum_{e \in E_P(\xi)} \norm{\pi_\xi(q_e - p_e)}_2} \\
        &= 2\pi\,\ExY[\big]{\xi}{\ExY{u \in \xi^\perp \cap \Sphere}{\hsupp{P}(u)}} \\
        &= 2\pi\,\ExY{u}{\hsupp{P}(u)} \\
        &= \pi\,\mw{P}.
    \end{align*}
    Multiplying both sides by $\frac{4}{\pi}$ yields the second identity.
\end{proof}

\begin{remark:unnumbered}
    The equality $\ExY{\xi}{\lenCY{2}{\Shad{P}{\xi}}} = 4\,\mw{P}$ in
    \lemref{polyhedral-crofton} is specific to polyhedral silhouettes: for
    a smooth body such as $B^3$, every shadow boundary has Euclidean length
    $2\pi < 8 = 4\,\mw{B^3}$ because polyhedral silhouette edges tilt
    relative to $\xi^\perp$, so only the lower-semicontinuity inequality
    passes to the limit in \thmref{short-shadow}.
\end{remark:unnumbered}

\begin{theorem}[Existence of a short shadow boundary]
    \thmlab{short-shadow}
    Suppose $B^3 \subseteq \CSym$ and the mean width satisfies $\mw{\CSym} < 3$.
    Then there exists a segment-free direction $\xi \in \SegFree{K}$
    such that $\lenCY{\CSym}{\Shad{K}{\xi}} < 12$.
\end{theorem}

\begin{figure}[!htb]
    \centering
    \includegraphics[width=0.88\textwidth]{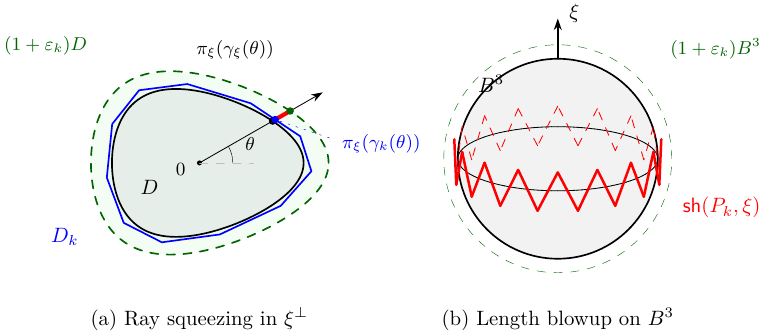}
    \caption{Length transfer in \thmref{short-shadow}: (a) ray squeezing of
        $\pi_\xi(\gamma_k(\theta))$ between $D$ and $(1+\varepsilon_k)D$; (b)
        equatorial zig-zag on $B^3$ showing shadow-boundary lengths need not
    converge.}
    \figlab{short-shadow-limit}
\end{figure}

\begin{proof}
    Set $\eta = 12 - 4\,\mw{\CSym} > 0$. For
    each integer $k \ge 1$, set $\varepsilon_k = \frac{1}{k}$ and apply
    \lemref{outer-polytope} to choose a convex polytope $P_k$ with
    $0 \in \interior(P_k)$ and $K \subseteq P_k \subseteq (1 +
    \varepsilon_k)K$. Since $B^3 \subseteq \CSym$, we have
    $\norm{v}_{\CSym} \le \norm{v}_2$ for all $v \in \R^3$; since $P_k
    \subseteq (1 + \varepsilon_k)K$, we have $\hsupp{P_k}(u) \le (1 +
    \varepsilon_k)\hsupp{K}(u)$ for all $u \in \Sphere$, and hence
    $\mw{P_k} \le (1 + \varepsilon_k)\mw{K} = (1 + \varepsilon_k)\mw{\CSym}$.
    By \thmrefY{segment-free-measure}{the theorem of Ewald \etal} and
    \obsref{polyhedral-silhouette},
    $\Omega = \SegFree{K} \cap \bigcap_{k=1}^\infty \GenDir{P_k}$ has
    full measure $\sigma(\Omega) = 1$. On $\Omega$,
    define $f_k(\xi) = \lenCY{\CSym}{\Shad{P_k}{\xi}}$, which is a
    nonnegative measurable step function by \obsref{polyhedral-silhouette}
    ($C = \CSym$). By \lemref{polyhedral-crofton} ($C = \CSym$) and the
    comparison inequalities above,
    \begin{align*}
        \ExY{\xi}{f_k(\xi)}
        &= \frac{1}{\pi}\sum_{e \in E(P_k)} \norm{q_e - p_e}_{\CSym}\theta(e)
        \le \frac{1}{\pi}\sum_{e \in E(P_k)} \norm{q_e - p_e}_2 \theta(e) \\
        &= 4\,\mw{P_k}
        \le 4(1 + \varepsilon_k)\,\mw{\CSym}
        = (1 + \varepsilon_k)(12 - \eta).
    \end{align*}
    By \lemrefY{fatou}{Fatou's lemma} (extending $f_k$ by $0$ outside the
    full-measure set $\Omega$),
    \begin{equation*}
        \ExY[\Big]{\xi}{\liminf_{k \to \infty} f_k(\xi)}
        \le \liminf_{k \to \infty} \ExY{\xi}{f_k(\xi)}
        \le \liminf_{k \to \infty} \pth{1 + \frac{1}{k}} (12 - \eta) =
        12 - \eta.
    \end{equation*}
    Since $\sigma(\Omega) = 1$, there exists $\xi \in \Omega \subseteq
    \SegFree{K}$ such that $\liminf_{k \to \infty} f_k(\xi) \le 12 - \eta
    < 12$.\footnote{The $\liminf$ is genuinely needed: the bounds
    $K \subseteq P_k \subseteq (1 + \varepsilon_k)K$ and
    $\xi \in \SegFree{K} \cap \bigcap_k \GenDir{P_k}$ do not imply that
    $f_k(\xi)$ converges to $\lenCY{\CSym}{\Shad{K}{\xi}}$ (or even stays
    bounded). Even for $K = B^3$ and the vertical direction
    $\xi = (0, 0, 1)$, taking the convex hull of a sufficiently tight outer
    polytope of $B^3$ and $2N$ alternating vertices at horizontal radius
    $1 + \varepsilon_k/4$ and height $\pm\sqrt{\varepsilon_k}$ around the
    equator yields a polytope $P_k \subseteq (1 + \varepsilon_k)B^3$ with
    $\xi \in \GenDir{P_k}$ whose shadow boundary zig-zags with length
    $> 4N\sqrt{\varepsilon_k}$, which can be arbitrarily large
    (\figref{short-shadow-limit}~(b)).}

    We are almost done---the closed curves $\Shad{P_k}{\xi}$ converge to
    $\Shad{K}{\xi}$ while their $\CSym$-lengths have limit inferior at most
    $12 - \eta$, so $\Shad{K}{\xi}$ also has $\CSym$-length at most
    $12 - \eta$. The legalese is as follows. Fix polar coordinates in
    $\xi^\perp$, let $D = \pi_\xi(K)$ and
    $D_k = \pi_\xi(P_k)$, and let $\gamma_\xi : [0, 2\pi] \to
    \Shad{K}{\xi}$ and $\gamma_k : [0, 2\pi] \to \Shad{P_k}{\xi}$ be the
    polar parameterizations from \obsref{outside-shadow}. Because $D
    \subseteq D_k \subseteq (1 + \varepsilon_k)D$
    (\figref{short-shadow-limit}~(a)), along each ray from
    $0$ at angle $\theta \in [0, 2\pi]$ we have $\pi_\xi(\gamma_k(\theta))
    \in [1, 1 + \varepsilon_k]\,\pi_\xi(\gamma_\xi(\theta))$, so
    $\pi_\xi(\gamma_k(\theta)) \to \pi_\xi(\gamma_\xi(\theta))$ as $k \to
    \infty$. For fixed $\theta \in [0, 2\pi]$, $\gamma_k(\theta) \in P_k
    \subseteq 2 K$ is bounded, and any cluster point $x_\infty$ of
    $\set{\gamma_k(\theta)}_{k \ge 1}$ satisfies $x_\infty = \lim_j (1 +
    \varepsilon_{k_j})^{-1}\gamma_{k_j}(\theta) \in K$ and
    $\pi_\xi(x_\infty) = \pi_\xi(\gamma_\xi(\theta))$, forcing $x_\infty =
    \gamma_\xi(\theta)$ by \obsref{outside-shadow-obs}. Thus
    $\gamma_k(\theta) \to \gamma_\xi(\theta)$ pointwise on $[0, 2\pi]$.

    For any partition $0 = \theta_0 < \theta_1 < \dots < \theta_N =
    2\pi$, continuity of $\norm{\cdot}_{\CSym}$ gives
    \begin{equation*}
        \sum_{r=1}^N \norm{\gamma_\xi(\theta_r) -
        \gamma_\xi(\theta_{r-1})}_{\CSym}
        = \lim_{k \to \infty} \sum_{r=1}^N \norm{\gamma_k(\theta_r) -
        \gamma_k(\theta_{r-1})}_{\CSym}
        \le \liminf_{k \to \infty} f_k(\xi) \le 12 - \eta.
    \end{equation*}
    Taking the supremum over all partitions yields
    $\lenCY{\CSym}{\Shad{K}{\xi}} = \lenCY{\CSym}{\gamma_\xi} \le 12 -
    \eta < 12$.
\end{proof}

\subsection{Putting everything together}
\seclab{main-proof}

\begin{proofof}
    \textit{\thmref{main}.} Let $K \subset \R^3$ be an arbitrary convex body. If $K$ is an
    affine cube, \obsref{parallelepipeds} and
    \propref{equiv-illumination} give $c(K) = 8$ with a uniform homothety
    ratio $\lambda \in (0, 1)$.

    Assume now that $K$ is not an affine cube. By \factref{affine-inv},
    after a translation and an invertible linear transformation we may
    assume that $0 \in \interior(K)$ and $\CSym = \frac{1}{2}(K - K)$ is in
    John's position ($B^3 \subseteq \CSym$). Since $K$ is not an affine
    cube, \thmref{mean-width-rigidity} gives $\mw{\CSym} < 3$.

    By \thmref{short-shadow}, there exists a direction $\xi \in
    \SegFree{K}$ with $L = \lenCY{\CSym}{\Shad{K}{\xi}} < 12$.
    Let $\gamma_\xi : [0, 2\pi] \to \Shad{K}{\xi}$ be its polar
    parameterization from \obsref{outside-shadow}. Since $\gamma_\xi$ is
    continuous with finite $\CSym$-length, the cumulative $\CSym$-length
    $V(\theta) = \lenCY{\CSym}{\gamma_\xi[0, \theta]}$ is continuous and
    nondecreasing on $[0, 2\pi]$ with $V(0) = 0$ and $V(2\pi) = L$.
    Applying the Intermediate Value Theorem gives parameters $0 = s_0 <
    s_1 < s_2 < s_3 < s_4 < s_5 < s_6 = 2\pi$ with $V(s_i) = \frac{i
    L}{6}$ for $i = 0, 1, \dots, 6$. For each $i \in [6]$, additivity
    gives
    \begin{equation*}
        \lenCY{\CSym}{\gamma_\xi[s_{i-1}, s_i]} = V(s_i) - V(s_{i-1}) =
        \frac{L}{6} < 2.
    \end{equation*}
    Since $\alpha_i = \gamma_\xi[s_{i-1}, s_i] \subseteq \Shad{K}{\xi}
    \subseteq \partial K$, \lemref{curve-covering} provides a nonzero
    vector $u_i \in \R^3 \setminus \set{0}$ that illuminates every point
    of $\alpha_i$ in $K$. Because $\bigcup_{i=1}^6 \alpha_i =
    \gamma_\xi[0, 2\pi] = \Shad{K}{\xi}$ by \obsref{outside-shadow}, the
    six vectors $u_1, \dots, u_6$ illuminate $\Shad{K}{\xi}$.

    Set $u_7 = +\xi$ and $u_8 = -\xi$. By \obsref{outside-shadow-obs},
    $\set{u_7, u_8}$ illuminates $\partial K \setminus \Shad{K}{\xi}$.
    Thus $u_1, \dots, u_8 \in \R^3 \setminus \set{0}$ illuminate all of
    $\partial K$. By \propref{equiv-illumination}, $c(K) \le 8$, and
    there exist $\lambda \in (0, 1)$ and $v_1, \dots, v_8 \in \R^3$ such
    that $K \subseteq \bigcup_{i=1}^8 (\lambda K + v_i)$.
\end{proofof}

\section{Concluding remarks}
\seclab{concluding}

A striking feature of the proof of \thmref{main} is that almost all of its
individual ingredients---Boltyanski's illumination equivalence, John's
decomposition, the mean-width bound and cube rigidity in John's position,
the polyhedral Crofton identity, the theorem of Ewald \etal, and Fatou's
lemma---are classical. What makes them yield the sharp bound $c(K) \le 8$
is that they fit together with no slack. First, \lemref{curve-covering}
covers any curve of $\CSym$-length strictly less than $2$ by a single
smaller homothet of $K$ itself, incurring no penalty for the asymmetry of
$K$ relative to its central symmetral $\CSym = \frac{1}{2}(K - K)$; the
entire metric and mean-width analysis can therefore take place in the
centrally symmetric geometry of $\CSym$. Second, in John's position for
$\CSym$, the norm comparison $\norm{\cdot}_{\CSym} \le \norm{\cdot}_2$,
the Crofton bound $4\,\mw{\CSym}$, and the mean-width ceiling
$4 \cdot 3 = 12$ match the six-arc illumination budget $6 \times 2 = 12$
exactly, with equality occurring only for affine cubes. Third,
transferring the length bound from outer polytopes to $\Shad{K}{\xi}$ via
Fatou's lemma and lower semicontinuity of curve length works directly on
$K$, bypassing the Hausdorff approximation trap for $c(K)$.

Our proof establishes $c(K) \le 8$ for all convex bodies $K \subset \R^3$,
leaving open the sharper question of whether every convex body in $\R^3$
other than an affine cube satisfies $c(K) \le 7$ (analogous to Levi's
characterization in the
plane). The averaging argument here naturally stops at $8$: for a slightly
smoothed cube, the central symmetral $\CSym$ in John's position has
$4\,\mw{\CSym}$ arbitrarily close to $12$, so our expectation bound can
only guarantee a shadow boundary of $\CSym$-length strictly less than
$12$---requiring six arcs of length less than $2$, even though a smooth
body satisfies $c(K) \le 4$.

Finally, the proof explains why three dimensions are special for this
approach. In any dimension $d \ge 3$, two opposite directions
$\set{+\xi, -\xi}$ illuminate $\partial K \setminus \Shad{K}{\xi}$,
leaving $2^d - 2$ directions for the shadow boundary $\Shad{K}{\xi}$. For
a segment-free direction $\xi$, the shadow boundary is homeomorphic to the
$(d - 2)$-dimensional sphere $S^{d-2}$. Only when $d = 3$ is this a
$1$-dimensional closed curve that can be partitioned by arc length and
covered piece by piece using \lemref{curve-covering}; already in $\R^4$,
the shadow boundary is a $2$-dimensional sphere, and controlling how many
smaller homothets of $K$ cover it requires a genuinely higher-dimensional
mechanism.

\bibliographystyle{alpha}%
\bibliography{refs}%

\appendix
\section{Deferred proofs}
\seclab{appendix}

\subsection{Proof of \texorpdfstring{\propref{equiv-illumination}}{equivalence of covering and illumination}}
\seclab{proof-equiv-illumination}
\begin{proof}
    The equivalence (i) $\iff$ (ii) is immediate (take $\lambda =
    \max_{i \in [m]} \lambda_i \in (0, 1)$ and enlarge each homothet around
    any point of $K$), as is (ii) $\implies$ (iii) (shrinking $K$ by
    $\lambda \in (0, 1)$ around any $p_0 \in \interior(K)$ places each
    $\lambda K + v_i$ inside $\interior(K) + w_i$ with $w_i = v_i - (1 -
    \lambda)p_0$).

    For (iii) $\implies$ (ii), fix $p_0 \in \interior(K)$ and set $K_0 =
    K - p_0$, so $0 \in \interior(K_0)$ and $K_0 \subseteq
    \bigcup_{i=1}^m (\interior(K_0) + w_i)$. By
    \lemref{gauge-props}, the function $f(y) = \min_{1 \le i
    \le m} \mu_{K_0}(y - w_i)$ is continuous on the compact set $K_0$ and
    satisfies $f(y) < 1$ for every $y \in K_0$. Hence $\lambda =
    \max\pth{\frac{1}{2}, \max_{y \in K_0} f(y)} \in (0, 1)$, and every
    $y \in K_0$ satisfies $\mu_{K_0}(y - w_i) \le \lambda$ for some $i
    \in [m]$, meaning $K_0 \subseteq \bigcup_{i=1}^m (\lambda K_0 +
    w_i)$. Adding $p_0$ gives $K \subseteq \bigcup_{i=1}^m (\lambda K +
    v_i)$ with $v_i = w_i + (1 - \lambda)p_0$.

    For (iii) $\implies$ (iv), set $u_i = -w_i$ if $w_i \neq 0$, and
    $u_i = e_1$ if $w_i = 0$. Given $x \in \partial K$, pick $i \in [m]$
    with $x \in \interior(K) + w_i$. Because $x \in \partial K = K
    \setminus \interior(K)$, we must have $w_i \neq 0$, so $u_i = -w_i
    \neq 0$ and $x + 1 \cdot u_i = x - w_i \in \interior(K)$.

    For (iv) $\implies$ (iii), for each $i \in [m]$ and $s > 0$ let
    $U_i(s) = \interior(K) - s u_i$. If $0 < s_1 < s_2$ and $x \in K
    \cap U_i(s_2)$, then $x \in K$ and $x + s_2 u_i \in \interior(K)$, so
    by convexity $x + s_1 u_i = (1 -
    \frac{s_1}{s_2})x + \frac{s_1}{s_2}(x + s_2 u_i) \in \interior(K)$
    and hence $x \in K \cap U_i(s_1)$. Thus the relatively open sets $W_n
    = \bigcup_{i=1}^m (K \cap U_i(1/n))$ ($n \ge 1$) satisfy $W_1
    \subseteq W_2 \subseteq \cdots$. Every $x \in \partial K$ satisfies
    $x + t_x u_i \in \interior(K)$ for some $i \in [m]$ and $t_x > 0$,
    placing $x \in W_n$ for all $n \ge 1/t_x$, while every $x \in
    \interior(K)$ has $B(x, \rho) \subseteq \interior(K)$ for some $\rho
    > 0$ and lies in $W_n$ for all $n > \norm{u_1}_2 / \rho$. By
    compactness of $K$, $K = W_N$ for some $N \ge 1$, giving $K \subseteq
    \bigcup_{i=1}^m (\interior(K) + w_i)$ with $w_i = -\frac{1}{N}u_i$.
\end{proof}

\subsection{Proof of \texorpdfstring{\thmref{john}}{John's decomposition}}
\seclab{proof-john}
\begin{proof}
    Let $C = \Sphere \cap \partial \CSym$ be the set of contact points, which is
    compact because $\Sphere$ and $\partial \CSym$ are compact. Moreover, $C \neq
    \emptyset$: otherwise $\Sphere \subset \interior(\CSym)$ (since $B^3 \subseteq
    \CSym$), so by \lemref{gauge-props} the continuous function $\norm{\cdot}_{\CSym}$
    satisfies $m_0 = \max_{u \in \Sphere} \norm{u}_{\CSym} < 1$, whence $m_0^{-1} B^3
    \subseteq \CSym$ with $m_0^{-1} > 1$, contradicting the maximality of
    $\vol_3(B^3)$ among origin-centered ellipsoids in $\CSym$.

    To verify~(ii), let $v \in C$. By \lemref{separation}~(ii), there exists
    $n \in \R^3$ with $\norm{n}_2 = 1$ such that $\iprod{n}{x} \le
    \iprod{n}{v}$ for all $x \in \CSym$. Since $n \in B^3 \subseteq \CSym$, we have
    \begin{equation*}
        1 = \iprod{n}{n} \le \iprod{n}{v} \le \norm{n}_2 \norm{v}_2 = 1,
    \end{equation*}
    so equality holds in Cauchy--Schwarz and hence $n = v$. Thus
    $\iprod{v}{x} \le 1$ for all $x \in \CSym$; replacing $x$ by $-x \in \CSym$
    (since $\CSym = -\CSym$) yields $\abs{\iprod{v}{x}} \le 1$ for all $x \in \CSym$.

    To establish~(i), equip the space $\Sym_3(\R)$ of real $3 \times 3$
    symmetric matrices with the Frobenius inner product $\iprod{A}{B}_F =
    \tr(AB)$ and norm $\norm{A}_F = \sqrt{\tr(A^2)}$. Since $C \subset \Sphere$
    is compact, the set $\Set{v v^\top}{v \in C} \subset \Sym_3(\R)$ is
    compact, and therefore its convex hull $\mathcal{H} = \conv\Set{v
    v^\top}{v \in C}$ is compact and convex.

    Suppose for contradiction that $\frac{1}{3} I_3 \notin \mathcal{H}$. By
    \lemref{separation}~(i) applied to $\mathcal{H} \subset \Sym_3(\R)$, there
    exists a nonzero matrix $W \in \Sym_3(\R)$ such that
    \begin{equation*}
        \tr(W X) \le \frac{1}{3}\tr(W) - \norm{W}_F^2
        \qquad\text{for all } X \in \mathcal{H}.
    \end{equation*}
    Set $\delta = \frac{1}{2}\norm{W}_F^2 > 0$ and define $H = W -
    \pth{\frac{1}{3}\tr(W) - \delta} I_3 \in \Sym_3(\R)$. Then $\tr(H) =
    3\delta > 0$ (so $\norm{H}_F > 0$), and for every $v \in C$,
    \begin{equation*}
        \iprod{v}{H v} = \tr(W v v^\top) - \frac{1}{3}\tr(W) + \delta \le
        -\delta < 0.
    \end{equation*}
    By continuity of $u \mapsto \iprod{u}{H u}$, the relatively open subset
    $U = \Set{u \in \Sphere}{\iprod{u}{H u} < -\delta/2}$ of $\Sphere$ contains
    $C$. Note also that $\norm{H u}_2 \le \norm{H}_F$ for every $u \in \Sphere$,
    as is immediate in an orthonormal eigenbasis of $H$.

    For every $u \in U$ and $0 < \varepsilon \le \delta / \norm{H}_F^2$, we
    have
    \begin{equation*}
        \norm{(I_3 + \varepsilon H)u}_2^2
        = 1 + 2\varepsilon\iprod{u}{H u} + \varepsilon^2\norm{H u}_2^2
        < 1 - \varepsilon\delta + \varepsilon^2\norm{H}_F^2 \le 1,
    \end{equation*}
    so $(I_3 + \varepsilon H)u \in B^3 \subseteq \CSym$. Meanwhile, on the
    compact set $F = \Sphere \setminus U \subseteq \interior(\CSym)$,
    \lemref{gauge-props} implies that $\norm{\cdot}_{\CSym}$ is continuous with
    $\norm{u}_{\CSym} < 1$ for all $u \in F$, and $B^3 \subseteq \CSym$ implies
    $\norm{x}_{\CSym} \le \norm{x}_2$ for all $x \in \R^3$. If $F \neq \emptyset$,
    set $m = \max_{u \in F} \norm{u}_{\CSym} < 1$; then for every $u \in F$ and
    $0 < \varepsilon \le (1 - m)/\norm{H}_F$,
    \begin{equation*}
        \norm{(I_3 + \varepsilon H)u}_{\CSym} \le \norm{u}_{\CSym} + \varepsilon\norm{H
        u}_{\CSym} \le m + \varepsilon\norm{H u}_2 \le m + \varepsilon\norm{H}_F
        \le 1,
    \end{equation*}
    so $(I_3 + \varepsilon H)u \in \CSym$. Thus, for sufficiently small
    $\varepsilon > 0$, we have $(I_3 + \varepsilon H)(\Sphere) \subseteq \CSym$, and
    by convexity $(I_3 + \varepsilon H)(B^3) \subseteq \CSym$. Writing
    $\lambda_1, \lambda_2, \lambda_3 \in \R$ for the eigenvalues of $H$, we
    have
    \begin{equation*}
        \det(I_3 + \varepsilon H)
        = \prod_{j=1}^3 (1 + \varepsilon\lambda_j)
        = 1 + (\lambda_1 + \lambda_2 + \lambda_3)\varepsilon + O(\varepsilon^2)
        = 1 + 3\delta\varepsilon + O(\varepsilon^2) > 1
    \end{equation*}
    for all sufficiently small $\varepsilon > 0$, contradicting the
    maximality of $\vol_3(B^3)$ in $\CSym$.

    Consequently, $\frac{1}{3} I_3 \in \mathcal{H}$. Since $\tr(v v^\top) =
    1$ for all $v \in C$, $\mathcal{H}$ lies in the $5$-dimensional affine
    hyperplane $\Set{X \in \Sym_3(\R)}{\tr(X) = 1}$. By Carath\'eodory's
    theorem, there exist $v_1, \dots, v_r \in C$ ($r \le 6$) and scalars
    $\alpha_1, \dots, \alpha_r > 0$ with $\sum_{s=1}^r \alpha_s = 1$ and
    $\sum_{s=1}^r \alpha_s v_s v_s^\top = \frac{1}{3} I_3$. Setting $c_s =
    3\alpha_s > 0$ gives $\sum_{s=1}^r c_s v_s v_s^\top = I_3$ and
    $\sum_{s=1}^r c_s = 3$. Because $(-v)(-v)^\top = v v^\top$, merging any
    antipodal pairs ensures $v_s \neq \pm v_t$ for $s \neq t$, and
    $\operatorname{rank}(I_3) = 3$ forces $r \ge 3$.
\end{proof}

\subsection{Proof of \texorpdfstring{\thmref{mean-width-rigidity}}{mean width in John's position}}
\seclab{proof-mean-width}
\begin{proof}
    For every $x \in \CSym$ and $u \in \R^3$, \thmref{john}~(i)--(ii) give
    \begin{equation*}
        \iprod{x}{u} = \sum_{s=1}^r c_s \iprod{v_s}{x}\iprod{v_s}{u}
        \le \sum_{s=1}^r c_s \abs{\iprod{v_s}{x}}\abs{\iprod{v_s}{u}}
        \le \sum_{s=1}^r c_s \abs{\iprod{v_s}{u}}.
    \end{equation*}
    Taking the maximum over $x \in \CSym$ yields $\hsupp{\CSym}(u) \le \sum_{s=1}^r c_s
    \abs{\iprod{v_s}{u}}$ for all $u \in \R^3$ (\figref{john-cube}). By
    $\Orth_3(\R)$-invariance of $\sigma$ and \Eqref{spherical-integral},
    every unit vector $v \in \Sphere$ satisfies $\ExY{u}{\abs{\iprod{v}{u}}} =
    \ExY{u}{\abs{u_3}} = \frac{1}{2}\int_0^\pi \abs{\cos\psi}\sin\psi \,
    d\psi = \frac{1}{2}$. Taking expectations over $u \in \Sphere$ therefore
    yields
    \begin{equation*}
        \mw{\CSym} = 2\,\ExY{u}{\hsupp{\CSym}(u)}
        \le 2 \sum_{s=1}^r c_s\,\ExY{u}{\abs{\iprod{v_s}{u}}}
        = \sum_{s=1}^r c_s = 3.
    \end{equation*}

    Suppose equality holds. Then $g(u) = \sum_{s=1}^r c_s
    \abs{\iprod{v_s}{u}} - \hsupp{\CSym}(u)$ is continuous and nonnegative on $\Sphere$
    with $\ExY{u}{g(u)} = 0$. If $g(u_0) > 0$ at some $u_0 \in
    \Sphere$, continuity gives $\delta \in (0, 2)$ and $\eta_0 > 0$ with $g
    \ge \eta_0$ on $\Sphere \cap \Bclosed(u_0, \delta)$, forcing $\ExY{u}{g(u)}
    \ge \eta_0 \frac{\delta^2}{4} > 0$ by
    \Eqref{cap-measure}. Thus $\hsupp{\CSym}(u) = \sum_{s=1}^r c_s
    \abs{\iprod{v_s}{u}}$ for all $u \in \Sphere$. Evaluating at $u = v_t$
    ($t \in [r]$) and using $\hsupp{\CSym}(v_t) = 1 = \norm{v_t}_2^2 = \sum_{s=1}^r
    c_s \iprod{v_s}{v_t}^2$ gives
    \begin{equation*}
        \sum_{s=1}^r c_s \abs{\iprod{v_s}{v_t}}\pth{1 -
        \abs{\iprod{v_s}{v_t}}} = 0.
    \end{equation*}
    Since $c_s > 0$ and $\abs{\iprod{v_s}{v_t}} \in [0, 1]$, we have
    $\abs{\iprod{v_s}{v_t}} \in \set{0, 1}$, and because $v_s \neq \pm
    v_t$ for $s \neq t$, we get $\iprod{v_s}{v_t} = 0$ for all $s \neq
    t$. Thus $v_1, \dots, v_r$ are orthonormal in $\R^3$, forcing $r = 3$
    and $c_t = \iprod{v_t}{I_3 v_t} = 1$. Defining $Q \in \Orth_3(\R)$
    with rows $v_1^\top, v_2^\top, v_3^\top$, the cube $Z = Q^\top([-1,
    1]^3) = \sum_{s=1}^3 [-1, 1]v_s$ satisfies $\hsupp{Z}(u) = \sum_{s=1}^3
    \abs{\iprod{v_s}{u}} = \hsupp{\CSym}(u)$ for all $u \in \R^3$. By
    \lemref{separation}~(i), two compact convex sets with identical support
    functions coincide, so $\CSym = Q^\top([-1, 1]^3)$ is an affine cube,
    and \lemref{diff-parallelepiped} shows $K$ is an affine cube.

    Conversely, suppose $K$ is an affine cube, so $\CSym = M([-1, 1]^3)$
    for some $M \in \GL_3(\R)$ by \lemref{diff-parallelepiped}. Since
    $B^3 \subseteq [-1, 1]^3$, we have $M(B^3) \subseteq \CSym$, and hence
    $\abs{\det M} \le 1$ by the maximality of $\vol_3(B^3)$ in $\CSym$. Set
    $T = M^{-1}$, with rows $r_1^\top, r_2^\top, r_3^\top$. Because
    $T(B^3) \subseteq [-1, 1]^3$, testing $u = r_i / \norm{r_i}_2 \in B^3$
    (if $r_i \neq 0$) gives $\norm{r_i}_2 \le 1$ for $i \in \set{1, 2,
    3}$. By Lagrange's cross-product identity $\norm{r_2 \times r_3}_2^2 =
    \norm{r_2}_2^2 \norm{r_3}_2^2 - \iprod{r_2}{r_3}^2$, we have
    \begin{equation*}
        1 \le \abs{\det T}^2 = \iprod{r_1}{r_2 \times r_3}^2 \le
        \norm{r_1}_2^2 \pth{\norm{r_2}_2^2 \norm{r_3}_2^2 -
        \iprod{r_2}{r_3}^2} \le \norm{r_1}_2^2 \norm{r_2}_2^2
        \norm{r_3}_2^2 \le 1,
    \end{equation*}
    forcing $\norm{r_i}_2 = 1$ and (applying the same equality to every
    row ordering) $\iprod{r_i}{r_j} = 0$ for $i \neq j$. Thus $T \in
    \Orth_3(\R)$, $\CSym = T^{-1}([-1, 1]^3) = T^\top([-1, 1]^3)$, and
    $\Orth_3(\R)$-invariance of $\sigma$ gives $\mw{\CSym} =
    2\,\ExY{u}{\abs{u_1} + \abs{u_2} + \abs{u_3}} = 3$.
\end{proof}

\end{document}